\documentclass[a4paper,12pt]{amsart}

\theoremstyle{definition}

\newtheorem{problem}{Problem}[section]

\usepackage{graphicx}

\usepackage{amsmath}
\usepackage{amssymb}
\usepackage{amsthm}
\usepackage{a4wide}
\usepackage{xcolor}

\usepackage{hyperref}
\usepackage{fancyhdr}
\usepackage{fourier}

\DeclareMathOperator{\lip12}{Lip_{1/2}}

\newcommand{\N}{{\mathbb N}}

\newcommand{\D}{{\mathbb D}}

\newcommand{\T}{{\mathbb T}}

\newcommand{\Ht}{{\mathcal{H}}^2}
\newcommand{\Hp}{{\mathcal{H}}^p}
\newcommand{\Ho}{{\mathcal{H}}^1}
\begin{document}

\author[Eero Saksman]{Eero Saksman}
\address{Department of Mathematics and Statistics \\
University of Helsinki \\
FI-00170 Helsinki \\ Finland} 
\email{eero.saksman@helsinki.fi}

\author[Kristian Seip]{Kristian Seip}
\address{Department of Mathematical Sciences \\ Norwegian University of Science and Technology \\ NO-7491 Trondheim \\ Norway}
\email{seip@math.ntnu.no}

\thanks{Saksman is
supported by the Finnish CoE in Analysis and Dynamics Research and
by a Knut and Alice Wallenberg Grant. Seip is supported by Grant 227768 of the Research Council of Norway.}

\title{Some open questions in analysis for Dirichlet series}


\maketitle

\begin{abstract} 
We present some open problems and describe briefly some possible research directions in the emerging theory of Hardy spaces of Dirichlet series and their intimate counterparts, Hardy spaces on the infinite-dimensional torus. Links to number theory are emphasized throughout the paper.
\end{abstract}

\vspace{5mm}

\section{Introduction}

We have in recent years seen a notable growth of interest in
certain functional analytic aspects of the theory of ordinary
Dirichlet series
\[ \sum_{n=1}^\infty a_n n^{-s}. \]
Contemporary research in this field owes much to the following
fundamental observation of H. Bohr \cite{BOHR}: By the
transformation $z_j=p_j^{-s}$ (here $p_j$ is the $j$th prime
number) and the fundamental theorem of arithmetic, an ordinary
Dirichlet series may be thought of as a function of infinitely many
complex variables $z_1, z_2, ...$. More precisely, in the Bohr correspondence,
 \begin{equation}\label{eq:bohr}
	F(s):= \sum_{n=1}^\infty a_{n} n^{-s}\;\sim\;  f(s):=\sum_{\nu \in \N_{\operatorname{fin}}^\infty} \widetilde a_\nu z^{\nu}, 
\end{equation}
where $n = p_1^{\nu_1} \cdots p_k^{\nu_k}$ and we identify $\widetilde a_\nu$ with the corresponding coefficient $a_n$,
and   $\N_{\operatorname{fin}}^\infty$ stands for the finite  sequences of positive indices.
 By a classical approximation
theorem of Kronecker, this is much more than just a formal
transformation: If, say, only a finite number of the coefficients
$a_n$ are nonzero (so that questions about convergence of the series
are avoided), the supremum of the Dirichlet polynomial $\sum a_n
n^{-s}$ in the half-plane $\operatorname{Re} s>0$ equals the
supremum of the corresponding polynomial on the infinite-dimensional
polydisc $\D^\infty$. In a groundbreaking work of Bohnenblust and
Hille \cite{BH}, it was later shown that homogeneous
polynomials---the basic building blocks of functions analytic on
polydiscs---may, via the method of polarization, be transformed into
symmetric multilinear forms. Bohnenblust and Hille used this insight to solve a long-standing problem in the field: 
Bohr had shown   that the width of the strip in which a Dirichlet  series converges uniformly but not absolutely is $ \leq 1/2$,
but  Bohnenblust and Hille   were able to prove that this upper estimate is in fact optimal.

In retrospect, one may in the work of Bohr and Bohnenblust--Hille
see the seeds of a theory of Hardy $H^p$ spaces of Dirichlet series.
However, this research took place before the modern interplay
between function theory and functional analysis, as well as the
advent of the field of several complex variables, and the area was
in many ways dormant until the late 1990s. One of the main goals of
the 1997 paper of Hedenmalm, Lindqvist, and Seip \cite{HLS97} was to
initiate a systematic study of Dirichlet series from the point of
view of modern operator-related function theory and harmonic
analysis. Independently, at the same time, a paper of Boas and
Khavinson \cite{BOKH} attracted renewed attention, in the context of
several complex variables, to the original work of Bohr.

The main object of study in \cite{HLS97} is the Hilbert space of
Dirichlet series $\sum_n a_n n^{-s}$ with square summable
coefficients $a_n$. This Hilbert space $\Ht$ consists of functions
analytic in the half-plane $\operatorname{Re} s>1/2$. Its
reproducing kernel at $s$ is $k_s(w)=\zeta(\overline{s}+w)$, where
$\zeta$ is the Riemann zeta function. By the Bohr correspondence, $\Ht$
may be thought of as the Hardy space $H^2$ on the
infinite-dimensional torus $\T^\infty$. Bayart \cite{Bayart02} extended the definition to any $p>0$ by defining $\Hp$ as the closure of Dirichlet polynomials $F(s)=\sum_{n=1}^{N}a_n n^{-s}$ under the norm
$$
\|F\|_{\Hp}\;:=\; \Big(\lim_{T\to\infty}\frac{1}{2T}\int_{-T}^T|F(it)|^p dt\Big)^{1/p}.
$$
By ergodicity (or see \cite{SaSe} for an elementary argument),  the Bohr correspondence yields the identity
\begin{equation}\label{eq:norm}
\|F\|_{\Hp}=\|f\|_{H^p(\T^\infty)}:=\Big(\int_{\T^\infty}|f(z)|^pd  m_\infty(z)\Big)^{1/p},
\end{equation}
where $ m_\infty$ stands for the Haar measure on the distinguished boundary $\T^\infty$, i.e., for the product of countably many copies of normalized Lebesgue measure on the circle $\T$.
Since the Hardy spaces on the infinite dimensional torus $H^p(\T^\infty)$  may be defined as the closure of analytic polynomials in the $L^p$-norm on $\T^\infty$, it follows that the Bohr correspondence provides an isomorphism between the spaces $H^p(\T^\infty)$ and $\Hp$. This linear isomorphism is  both isometric and multiplicative.

The classical theory of Hardy spaces and the operators that act on
them serves as an important source of incitement for the field of
Dirichlet series that has evolved after 1997. Two distinct features
should however be noted. First, a number of new phenomena, typically
crossing existing disciplines, appear that are not present in the
classical situation. Second, many of the classical objects change
radically and require new viewpoints and methods in order to be
properly understood and analyzed. 

In the following sections, we  sketch briefly some 
research directions and list several open problems (thus updating \cite{Hed}). In our selection of problems, we have followed our own interests and made no effort to compile a comprehensive list. As a consequence, several interesting recent developments such as for instance \cite{BDFMS} or \cite{MS} will not be accounted for and discussed. The reader should also notice that the difficulty of the problems may vary considerably. 
It seems likely that for some of the problems mentioned below,
further progress will require novel and unconventional combinations of tools from harmonic, functional, and complex analysis, as well as from
analytic number theory.

\section{Basic properties of the spaces $\Hp$ and $H^p(\T^\infty)$} \label{se:basic}

The study of the boundary limit functions in the spaces $\Hp$ has a
number of interesting features. Several central points have been
clarified, such as questions concerning convergence of the Dirichlet
series \cite{HS}, to what extent ergodicity extends to the boundary
\cite{SaSe}, properties of the boundary limit functions for
Dirichlet series in $\Ht$ \cite{OlSa}, and zeros of functions in $\Ht$ and, at least partially, in $\Hp$ for $p> 2$ \cite{Se}. The diversity of techniques
involved is considerable, ranging from function theory in polydiscs
and ergodic theory to classical harmonic analysis, Hardy space
techniques, Fourier frames, estimates for solutions of the $\overline{\partial}$ equation, and Ramanujan's estimates for the divisor function. Still, a very natural problem first considered
in \cite{Bayart02}  (see \cite{SaSe} for further discussion on it) remains unsolved and represents
one of the main obstacles to further progress: \vspace{1mm}

\begin{problem}[The embedding problem]\label{embedding}
Is the $L^p$ integral of a Dirichlet
polynomial $\sum_{n=1}^N a_n n^{-s}$ over any segment of fixed length on the
vertical line $\operatorname{Re} s=1/2$ bounded by a universal
constant times
$\| \sum_{n=1}^N a_n n^{-s} \|_{\Hp}^{p}$?
\end{problem}

\noindent This is known to hold for $p=2$ and thus trivially for $p$ an even
integer. One may notice a curious resemblance with Montgomery's conjectures concerning norm
inequalities for Dirichlet polynomials (see \cite[pp. 129, 146]{Mo94} or \cite[p. 232--235]{IK05}). It remains to be clarified if
there is a link between this question and Montgomery's conjectures. 

An affirmative answer to Problem \ref{embedding}  for $p<2$ would have immediate function theoretic consequences regarding for instance zero sets and boundary limits. Namely,  following \cite{OlSa}, we would be able to answer

\begin{problem}\label{Carleson}
Characterize Carleson measures for $\Hp$ on $\{\operatorname{Re} s>1/2\}$ for $p<2.$
\end{problem}

More modest but nontrivial open questions are:
 
 \begin{problem}\label{blaschke}
Do the zero sets of functions in $\Hp$ for $p<2$ satisfy the Blaschke condition in the half-plane $\operatorname{Re} s>1/2$?  
\end{problem}
 \begin{problem}\label{nevanlinna}
Are elements of $\Hp$ for $p<2$ locally in the Nevanlinna class?
\end{problem}

There are similar problems of a dual flavor regarding interpolating sequences for $\Hp$. Indeed, it follows from \cite{OlSe} that the Shapiro--Shields version of Carleson's classical theorem in the half-plane $\operatorname{Re} s>1/2$ remains valid when $1/p$ is an even integer. We would like  to know if this result extends to other values of $p$.

By a theorem of Helson \cite{H2}, the partial sum operator \cite{H2} is uniformly  bounded on  $\Hp$ for $1<p<\infty$ (see \cite{AOS} for an alternative treatment), and hence the functions $n^{-s}$ for $n\ge 1$ form a basis for $\Hp$ for these exponents $p$.  The following questions stated in \cite{AOS} seem to be open:

 \begin{problem}\label{unconditional}
Does $\Hp$ have an unconditional basis if $p \in (1,\infty)$ and $p\not=2$?
\end{problem}

 \begin{problem}\label{basis}
Does $\Ho$ have a basis? Does it have an unconditional basis?
\end{problem}

The last two problems are equivalent to corresponding statements for $H^p(\T^\infty).$ There are also natural and interesting questions that are  specific for function theory in infinite dimensions. In \cite{AOS2} (see \cite{SaSe} for the first steps in this direction), it was shown that Fatou or Marcin- kiewicz--Zygmund-type theorems on boundary limits remain true for all  classes $H^p(\T^\infty)$ or for their harmonic counterparts $h^p(\T^\infty)$, assuming fairly regular radial approach to the distinguished boundary $\T^\infty$; the simplest example of such approach is of the form $(re^{i\theta_1},r^2e^{i\theta_2},$ $r^3e^{i\theta_3},\ldots )$ with $r\uparrow 1^{-}.$ However, \cite{AOS2} also constructs an example of an element $f$ in $ H^\infty (\T^\infty)$ such that at almost every boundary point, $f$ fails to have  a radial limit under a certain radial approach that is independent of the boundary point.

 \begin{problem}\label{fatou}
Give general conditions for a radial (or non-tangential) approach in $\D^\infty$ to $\T^\infty$ such that Fatou's theorem holds for elements in $H^p(\T^\infty)$.
\end{problem}

The $\Hp$ spaces are well defined (via density of polynomials) also in the range $0<p<1$. Again, one may inquire the analogue of the embedding problem (now stated in term of local Hardy spaces on $\operatorname{Re} s=1/2$). For all values other than $p=2$, even partial non-trivial results pertaining to the following widely open question (see \cite{SaSe}) would be interesting.

 \begin{problem}\label{duals}
Describe the dual spaces of $\Hp$. 
\end{problem}


 \section{Operator theory and harmonic analysis} 
Viewing our Hardy spaces as closed subspaces of the ambient $L^p$
spaces on the infinite-dimensional torus $\T^\infty$, we are led to
consider classical operators like the Riesz projection (orthogonal
projection from $L^2$ to $H^2$),  Hankel operators, and Fourier multiplier operators. \cite{AOS} contains some results on multipliers
and Littlewood--Paley decompositions. 
It has become clear, however, that most of the classical methods are either
not relevant or at least insufficient for the infinite-dimensional
situation. For example,
the classical Nehari theorem for Hankel forms (or small
Hankel operators) does not carry over to $\T^\infty$, see \cite{OrSe}.
This leads us to ask if a reasonable replacement can be found and, more generally, how the different roles and interpretations of BMO (the space of functions of bounded mean oscillation) manifest themselves in our infinite-dimensional setting.

\begin{problem}\label{nehari} What is the counterpart to
Nehari's theorem on $\T^\infty$? In particular, what can be said about the Riesz projection of $L^\infty(\T^\infty)$ and other BMO-type spaces on $\T^\infty$?
\end{problem}

This and similar operator theoretic problems may be approached along several different paths. In \cite{BPSSV},  a natural analogue of the classical Hilbert matrix was identified and studied. This matrix was referred to as the multiplicative Hilbert matrix because its entries $a_{m,n}:=(\sqrt{mn} \log(mn))^{-1}$ depend on the product $m\cdot n$. This matrix represents a bounded Hankel form on $\Ht_0 \times\Ht_0 $ with spectral problems similar to those of the classical Hilbert matrix. (Here $\Ht_0$ denotes the subspace of $\Ht$ consisting of functions that vanish at $+\infty$.) Its analytic symbol $\varphi_0$ is a primitive of $-\zeta(s+1/2)+1$, and by analogy with the classical situation, we are led to the following problem.

\begin{problem}\label{nehari2} 
Is the symbol $\varphi_0(s)=1+\sum_{n=2}^\infty (\log n)^{-1}n^{-1/2-s}$   the Riesz projection of a function in $L^\infty(\T^\infty)$?
\end{problem} 

It is interesting to notice that a positive answer to Problem \ref{embedding} for $p=1$ would yield a positive answer to this question, via an argument involving Carleson measures. We refer to \cite{BPSSV} for details.

The beautiful pioneering contribution of Gordon and Hedenmalm \cite{GH} and a growing number of other papers have established the study of composition operators on Hardy spaces of Dirichlet series as an active research area in the interface of one and several complex variables. In the series of papers \cite{QS, BQS, BB}, quantitative and functional analytic tools have been developed in this context, for example norm estimates for linear combinations of reproducing kernels, Littlewood--Paley formulas, and (soft) functional analytic remedies for the fact that $\Hp$ fails to be complemented when $1\le p < \infty$ and $p\neq 2$. 

\begin{problem}\label{composition} 
Characterize the compact composition operators on $\Ht$.
\end{problem} 

\section{Moments of sums of random multiplicative functions}

There has during the last few years been an interesting interplay between the study of sums of random multiplicative functions and problems and methods coming from Hardy spaces. This topic has a long history, beginning with an import paper of Wintner \cite{W}. One of the links to Hardy spaces comes from 

\begin{problem}[Helson's problem \cite{He}]\label{helson} Is it true that $\| \sum_{n=1}^N n^{-s}\|_{\Ho}= o(\sqrt{N})$ when $N\to\infty$. 
\end{problem}
This intriguing open problem arose from Helson's study of Hankel forms and a comparison with the one-dimensional Dirichlet kernel. However, it seems to be more fruitful to think of the problem in probabilistic terms, viewing the functions $p_j^{-s}$ as independent Steinhaus variables.  Resorting to a decomposition into homogeneous polynomials and using well known estimates for the arithmetic function $\Omega(n)$, it was shown in \cite{BS2} that 
$\| \sum_{n=1}^N n^{-s} \|_{\Ho} \gg \sqrt{N} (\log N)^{-0.05616} $. This was later improved by Harper, Nikeghbali, and Radziwi{\l\l} \cite{HNR} who, using methods from \cite{H}, found the lower bound $\sqrt{N} (\log\log N)^{-3+o(1)}$. In a recent preprint \cite{HL}, Heap and Lindqvist made a prediction based on  random matrix theory that Helson's conjecture is false. 

The preprint \cite{BS2} also gave a precise answer to the question of for which $m$ the homogeneous 
Dirichlet polynomials $\sum_{\Omega(n)=m, n\le N} n^{-s}$  have comparable $L^4$ and $L^2$ norms. Indeed, this happens if and only $m$ is, in a precise sense, strictly smaller that $\frac{1}{2}\log\log N$. An interesting problem coming from analytic number theory and the work of Hough \cite{Ho}, is to extend this result to higher moments. 

\begin{problem}\label{hough} Assume $k$ is an integer larger than $1$. For which $m$ (depending on $N$) will the $L^{2k}$ norms of $m$-homogeneous Dirichlet polynomials of length $N$ be comparable to their $L^2$ norms?
\end{problem}

Cancellations in the partial sums of the Riemann zeta function on the critical line can be studied through a similar problem concerning $\Hp$ norms. 
\begin{problem}\label{critical _polynomials} Determine the asymptotic behavior of $\big\| \sum_{n=1}^N n^{-1/2-s}\big\|_{\Hp}$ when $N\to \infty$ for  $0<p\leq 1.$
\end{problem}
An interesting modification of this problem is the following.
\begin{problem}\label{critical _polynomials2} Determine the precise asymptotic growth of  $\| \sum_{n=1}^N [d(n)]^\gamma n^{-1/2-s} \|_{\Hp}$ when $N\to \infty$  for $p\le 1$. 
\end{problem}
\noindent  A more general problem is to do the same for polynomials with coefficients represented by  multiplicative functions satisfying appropriate growth conditions.\cite{BHS} established the inequality 
\begin{equation} \label{below} \left(\sum_{n=1}^{N} |\mu(n)| |a_n|^2 [d(n)]^{\frac{\log p}{\log 2}-1}\right)^{1/2}\le \| f\|_{\Hp}, \end{equation}
valid for $f(s)=\sum_{n=1}^N a_n n^{-s}$ and $0<p\le 2$, where $\mu(n)$ is the M\"{o}bius function. This inequality, which should be recognized as an $L^p$-analogue of an inequality of Helson \cite{H2}, yields the lower bound
\begin{equation}\label{euler} \Big\|  \sum_{n=1}^N n^{-1/2-s}\Big\|_{\Hp} \gg (\log N)^{p/4} \end{equation}
for all $0<p<\infty$. An estimate in the opposite  direction in the range $1<p<\infty$ follows by applying Helson's theorem on the $L^p$ boundedness of the partial sum operator \cite{H2} on suitably truncated Euler product. When $p=1$ the same method yields that an additional factor $\log\log N$ appears on the right-hand side when $\gg$ is replaced by $\ll$ in \eqref{euler}, and thus Problems \ref{critical _polynomials} and \ref{critical _polynomials2} remain open exactly in the range $p\leq 1.$ Some results for Problem \ref{critical _polynomials2} are contained in the manuscript \cite{BSS}.

 A closely related and more general problem concerns the natural partial sum operator of the Dirichlet series whose $L^p$ norm can be estimated by Helson's theorem \cite{H2} for finite $p$ and a result from \cite{BCQ06} for $p=\infty$.

\begin{problem}\label{partialsumoperator} Determine the precise asymptotic growth of the norm of the partial sum operator
 $S_N:\sum_{n=1}^\infty a_n n^{-s}\mapsto 
\sum_{n=1}^N a_n n^{-s} $ when $N\to \infty$  for  $p=1$ (or more generally, for $p\le 1$ or $p=\infty$). 
\end{problem}
In the case $p=1$, a trivial one dimensional estimate yields a lower bound of order $\log\log N$, whereas \cite{BSS}  gives an upper bound of order $\log N/\log\log N$, so that presently  there is a large gap between the  known bounds.

We finish this section by recalling a pointwise version of the analogue of Helson's problem on the torus. Thus, for primes $p$ let $\chi(p):$ we i.i.d random variables with uniform distribution on $\T$ and define $\chi(n)=\prod_{k=1}^\ell \chi(p_k)^{\ell_l}$ for $n=p_1^{k_1}\ldots p_\ell^{k_\ell}.$

\begin{problem}\label{erdos_characters} Determine the almost sure growth rate  (in $N$) of the character sum
$$
\sum_{n=1}^N\chi(n).
$$
\end{problem}
This problem stems from Wintner, and is listed by Erd\H{o}s, although in the original version instead $\chi(p)$:s are Rademacher variables. Deep results on the problem were provided by Halasz \cite {Ha} in the 1980s, and recently Harper \cite{H} obtained remarkable improvements for the lower bound. But the original problem remains open.

 \section{Estimates for GCD sums and the Riemann zeta function}
 
The study of greatest common divisor (GCD) sums of the form
\begin{equation}\label{gcda}
\sum_{k,\ell=1}^N\frac{(\gcd(n_k,n_{\ell}))^{2\alpha}}{(n_k
n_{\ell})^\alpha}
\end{equation}
for $\alpha>0$ was initiated by Erd\H{o}s who inspired G\'{a}l \cite{G} to solve 
a prize problem of the Wiskundig Genootschap in Amsterdam concerning the case $\alpha=1$. 
G\'{a}l proved that when $\alpha=1$, the optimal upper bound for \eqref{gcda} is $CN(\log \log N)^2$, with $C$ an absolute constant independent of $N$ and the distinct positive integers $n_1,...,n_N$. The problem solved by G\'{a}l had been posed by Koksma in the 1930s, based on the observation that such bounds would have implications for the uniform distribution of sequences $(n_k x)$ mod 1 for almost all $x$.

Using the several complex variables perspective of Bohr and seeds found in \cite{LiSe}, Aistleitner, Berkes and Seip \cite{ABS} proved sharp upper bounds for \eqref{gcda} in the range $1/2<\alpha <1$ and a much improved estimate for $\alpha=1/2$, solving in particular a problem of Dyer and Harman \cite{DH}. The method of proof was based on identifying \eqref{gcda} as a certain Poisson integral on $\D^\infty$. The acquired bounds were also used to establish a Carleson--Hunt-type inequality for systems of dilated functions of bounded variation or belonging to $\lip12$, a result that in turn settled two longstanding problems on the almost everywhere behavior of systems of dilated functions. The Carleson--Hunt inequality and the original inequality of G\'{a}l (see \eqref{aone} below) were later optimized by Lewko and Radziwi{\l\l} \cite{LR}.  

Additional techniques were introduced by Bondarenko and Seip \cite{BS1, BS2} to deal with the limiting case $\alpha=1/2$, and finally the range $0<\alpha<1/2$ was clarified in \cite{BHS}. Writing
\[ \Gamma_\alpha(N):=
\frac{1}{N} \sup_{1\le n_1<n_2<\cdots < n_N} \sum_{k,\ell=1}^N\frac{(\gcd(n_k,n_{\ell}))^{2\alpha}}{(n_k
n_{\ell})^\alpha}, \]
we may summarize the state of affairs as follows:
\begin{align} \Gamma_1(N) & \sim \frac{6e^{2\gamma}}{\pi^2} \log\log N  \label{aone}\\
\log \Gamma_{\alpha}(N) & \asymp_{\alpha}  \frac{(\log N)^{(1-\alpha)}}{(\log\log N)^{\alpha}}, \quad 1/2<\alpha < 1\nonumber  \\
\log \Gamma_{1/2}(N) & \asymp  \sqrt{\frac{\log N \log\log\log N}{\log\log N}} \nonumber \\
\log \Gamma_{\alpha}(N)-(1-2\alpha)\log N & \asymp_{\alpha} \log\log N,  \quad 0<\alpha<1/2, \nonumber
\end{align}
where in \eqref{aone}, $\gamma$ denotes Euler's constant; these estimates remain the same if we replace $\Gamma_\alpha(N)$ by the possibly larger quantity  
\[ \Lambda_\alpha(N):=\sup_{1\le n_1<n_2<\cdots < n_N, \| c \|=1} \sum_{k,\ell=1}^Nc_k c_\ell \frac{(\gcd(n_k,n_{\ell}))^{2\alpha}}{(n_k
n_{\ell})^\alpha}, \]
where the vector $c=(c_1, c_2,..., c_N)$ consists of nonnegative numbers and $\| c\|^2:=c_1^2+c_2^2+\cdots c_N^2$.

Aistleitner  \cite{A} made the important observation that such estimates can be used to obtain $\Omega$-results for the Riemann zeta function. Indeed, using Hilberdink's version of the resonance method \cite{H}, he found a new proof of Montgomery's $\Omega$-results for $\zeta(\alpha+it)$ in the range $1/2<\alpha<1$ \cite{M}. In turn, Bondarenko and Seip applied the particular set $\{n_1, n_2,...,n_N\}$  yielding the lower bound for $\Lambda_{1/2}(N)$  in combination with the resonance method of Soundararajan \cite{So} 
 to obtain (unconditionally) the following:  given $c<1/\sqrt{2}$, there exists a $\beta$, $0<\beta<1$, such that for every sufficiently large $T$ 
 \begin{equation} \label{bonseip}
 \sup_{t\in (T^\beta , T)} |\zeta(1/2+it)|\ge \exp\left(c\sqrt{\frac{\log T\log\log\log T}{\log\log T}}\right). 
\end{equation}
This gives an improvement by a power of  $\sqrt{\log\log\log T}$ compared with previously known estimates \cite{BR,So}.  

We list 
two rather general questions pertaining to these recent developments.


\begin{problem} \label{bmoa}
Link the estimates for GCD sums to the function and operator theory of the spaces $\Hp$.
\end{problem}

\begin{problem} \label{moreriemann}
Develop further applications to and links with the Riemann zeta function. 
\end{problem}


Problem \ref{bmoa} originates in the observation from \cite{ABS} that GCD sums can be interpreted as Poisson integrals on polydiscs. Taking into account the prominent role played by Poisson integrals and the Poisson kernel in the classical setting (for instance in connection with functions of bounded mean oscillation), we are led to ask for potential function and operator theoretic interpretations or applications of our estimates for GCD sums.

Finally, we would like to give an example related to the rather vague and general Problem~\ref{moreriemann}. It concerns estimates relating the size of the coefficients to the $\Hp$ norm of a Dirichlet series, which can be traced back to Bohr's problem of computing the maximal distance between the abscissas of absolute and uniform convergence. Bohnenblust and Hille's solution to this problem  \cite{BH} relied on a revolutionary method of polarization for estimating the size of the coefficients of homogeneous polynomials. There was a revival of interest  in Bohnenblust and Hille's work after the 1997 paper of Boas and Khavinson \cite{BOKH} on so-called Bohr inequalities. It was gradually recognized that the original estimate of order $m^m$ for the constant $C(m)$ in the Bohnenblust--Hille inequality was not sufficiently accurate to reach the desired level of precision in various applications. Based on a re-examination of the original proof, a sophisticated version of H\"{o}lder's inequality due to Blei \cite{Bl}, and a Khinchin-type inequality of Bayart \cite{Bayart02}, Defant, \mbox{Frerick}, Ouna\"{\i}es,
    Ortega-Cerd\`{a}, and Seip established in  \cite{DFOOS11} that $C(m)$ grows at most exponentially in $m$. This was recently improved further by Bayart, Pellegrino, and Seoane-Sep\'{u}lveda \cite{BPS} who were able to show, by taking a new approach to Blei's inequality, that $C(m)$ grows at most as $\exp\big(c\sqrt{m\log m}\big)$ for some constant $c$. 

The most important application of the improved version of the Bohnenblust--Hille inequality was to the compute the Sidon constant $S(N)$ which is defined as the supremum of the ratio between $|a_1|+\cdots+
|a_N|$ and $\sup_{t\in {\Bbb R}}\big|a_1+a_2 2^{it}+\cdots + a_N
N^{it}\big|$, with the supremum taken over all possible choices of
nonzero vectors $(a_1, ..., a_N)$ in ${\Bbb C}^N$. The following
remarkably precise asymptotic result holds \cite{DFOOS11}:
\[ S(N)=\sqrt{N} \exp\left(\big(-\frac{1}{\sqrt{2}}+o(1)\big)\sqrt{\log N \log \log N}\right) \]
when $N\to \infty$. This formula has a long history and relies on
the contribution of many researchers, most notably Queff\'{e}lec and
Konyagin \cite{KQ01} and de la Bret\`{e}che \cite{Br08}. The proof
involves an unconventional blend of techniques from
function theory on polydiscs (the Bohnen- blust--Hille inequality),
analytic number theory (the Dickman function), and probability (the Salem--Zygmund inequality). 

There is a striking resemblance between the formula  for the Sidon constant $S(N)$ and the following conjecture from \cite{FGH}, based on  arguments from random matrix theory, conjectures for moments of $L$-functions, and also by assuming a random model for the primes  \cite{FGH}:
\[ \max_{0\le t\le T} \big|\zeta(1/2+it)\big| =\exp\left(\big(-\frac{1}{\sqrt{2}}+o(1)\big)\sqrt{\log T \log \log T}\right)  \]
when $T\to \infty$. It is natural to ask if this resemblance is more than just a coincidence.

\section {Random Dirichlet series}

 A classical result due to Selberg states that the distribution of the Riemann zeta function on the critical line is asymptotically Gaussian, after suitable renormalisation. More precisely, the distribution of
$
\left\{ \left(\frac{1}{2}\log\log (T)\right)^{-1/2}\log|\zeta(1/2+it)|\; :\; t\in [0,T]\right\}
$
tends to that of a standard normal variable $\mathcal N(0,1)$ as $T\to\infty.$ Recently, Fyodorov, Keating and Hiary computed heuristically the covariance of the translations of the zeta function and observed that in the first approximation a logarithmic correlation structure emerges. Similar covariance structure is exhibited by (the one-dimensional) restriction of the Gaussian free field (GFF), a fundamental probabilistic object that figures prominently in e.g. Liouville quantum gravity, SLE and random matrix theory. Based on the classical (after Montgomery) heuristic connection between $\zeta (1/2+it)$ and random matrices, and the conjectured behaviour of random matrices they proposed the following \begin{problem}\label{FHK}\cite{FHK}
Consider $[0,T]$ as a probability space, with normalised Lebesgue measure, and denote the corresponding variable 
by $\omega\in [0,T]$. Then, as $T\to\infty,$ one has
$$
\max_{h\in [0,1]}\log |\zeta (1/2+ih +i\omega) |=\log\log T -\frac{3}{4}\log\log\log T + E,
$$
where the error term $E$ is bounded in probability as $T\to\infty.$
\end{problem}
\noindent Very recently Arguin, Belius and Harper \cite{ABH} established the analogue of the above conjecture for a natural model that is derived  from the Euler product of the zeta-function, i.e. for partial sums of  random Dirichlet series  of the type
\begin{equation}\label{toy}
\notag X(x)=\sum_{p}\frac{1}{\sqrt{p}}\left(\cos(x\log p)\cos \theta_p+\sin(x\log p)\sin\theta_p\right),
\end{equation}
where $\theta_p:s $ are i.i.d. and unifrom on  $[0,2\pi]$ and indexed by prime numbers.

The GFF heuristics of the zeta-function over the critical line has been useful also in connection with the Helson conjecture \cite{HNR}. Many fascinating questions remain to be studied in this general domain of probabilistic behaviour of the zeta function and related models.  For many random Gaussian fields (taking values in generalised functions) one may construct the corresponding multiplicative Gaussian chaos measure
see e.g. \cite{K}, \cite{DRSV}, \cite{BKNSW}). Naturally, after Selberg's result one may inquire if one could produce a gaussian chaos as a suitable scaling limit of the Riemann zeta function on the critical line. An easier task would be to consider
\begin{problem}\label{FHK}\cite{FHK}
Let the field $X$ be defined as in \eqref{toy}. Study the properties non-Gaussian chaos "$\exp (\beta X(x))$".
\end{problem}
\noindent 
Some very early steps in this direction are contained in \cite{SW}.


\end{document}